\documentclass[12pt]{article}
\usepackage{amssymb}
\usepackage{amsmath,amsthm}
\usepackage[latin1]{inputenc}
\usepackage{hyperref}
\usepackage{color}
\usepackage{graphicx}
\hypersetup{colorlinks=true, linkcolor=blue, citecolor=blue,
urlcolor=blue}

 \newtheorem{remark}{Remark}

 \newtheorem{lemma}[remark]{Lemma}
 \newtheorem{theorem}[remark]{Theorem}
 
 \newtheorem{corollary}[remark]{Corollary}
  \newtheorem{claim}[remark]{Claim}

\title{On the metric dimension of corona product graphs}

\author{I. G.
Yero$^{1}$, D. Kuziak$^{2}$ and  J. A.
Rodr\'{\i}guez-Vel\'{a}zquez$^{1}$\\
    \\
$^1${\small Departament d'Enginyeria Inform\`atica i Matem\`atiques,}\\
{\small Universitat Rovira i Virgili,}  {\small Av. Pa\"{\i}sos
Catalans 26, 43007 Tarragona, Spain.} \\{\small
ismael.gonzalez\@@urv.cat, juanalberto.rodriguez\@@urv.cat}
\\
$^2${\small Faculty of Applied Physics and  Mathematics}\\
{\small Gda\'nsk University of Technology,} {\small
ul. Narutowicza 11/12 80-233 Gda\'nsk, Poland} \\ {\small
dkuziak\@@mif.pg.gda.pl}\\
}

\begin{document}

\maketitle

\begin{abstract}
Given a set of vertices $S=\{v_1,v_2,...,v_k\}$ of a connected graph $G$, the metric representation of a vertex $v$ of $G$ with respect to $S$ is the vector $r(v|S)=(d(v,v_1),d(v,v_2),...,d(v,v_k))$, where $d(v,v_i)$, $i\in \{1,...,k\}$ denotes the distance between  $v$ and $v_i$. $S$ is a resolving set for $G$ if for every pair of vertices $u,v$ of $G$, $r(u|S)\ne r(v|S)$. The metric dimension of $G$,  $dim(G)$, is the minimum cardinality of any resolving set for $G$. Let $G$ and $H$ be two graphs of order $n_1$ and $n_2$, respectively. The corona product $G\odot H$ is defined as the graph obtained from $G$ and $H$ by taking one copy of $G$ and $n_1$ copies of $H$ and  joining by an edge each vertex from the $i^{th}$-copy of $H$ with the $i^{th}$-vertex of $G$. For any integer  $k\ge 2$, we define the graph $G\odot^k H$ recursively from $G\odot H$ as $G\odot^k H=(G\odot^{k-1} H)\odot H$.
We give several results on the metric dimension of $G\odot^k H$. For instance, we show that given two connected graphs $G$ and $H$ of order $n_1\ge 2$ and $n_2\ge 2$, respectively, if the diameter of $H$ is at most two, then $dim(G\odot^k H)=n_1(n_2+1)^{k-1}dim(H)$. Moreover, if $n_2\ge 7$ and the diameter of $H$ is greater than five or $H$ is a cycle graph, then
$dim(G\odot^k H)=n_1(n_2+1)^{k-1}dim(K_1\odot H).$
\end{abstract}



\section{Introduction}


 The concepts of resolvability and location in graphs were described
independently by Harary and Melter \cite{harary} and Slater
\cite{leaves-trees}, to define the same structure in a
graph. After these papers were published several authors
developed diverse theoretical works about this topic
\cite{pelayo1,pelayo2,chappell,chartrand,chartrand1,landmarks,survey,tomescu}.
 Slater described the usefulness of these ideas into long range
aids to navigation \cite{leaves-trees}. Also, these concepts  have
some applications in chemistry for representing chemical compounds
\cite{pharmacy1,pharmacy2} or to problems of pattern recognition and
image processing, some of which involve the use of hierarchical data
structures \cite{Tomescu1}. Other applications of this concept to
navigation of robots in networks and other areas appear in
\cite{chartrand,robots,landmarks}. Some variations on resolvability
or location have been appearing in the literature, like those about
conditional resolvability \cite{survey}, locating domination
\cite{haynes}, resolving domination \cite{brigham} and resolving
partitions \cite{chappell,chartrand2,fehr,yerocartpartres}. In this article we study the metric dimension of corona product graphs.

We begin by giving some basic concepts and notations. Let $G=(V,E)$ be a simple graph of order $n=|V|$. Let $u,v\in V$ be two different vertices in $G$, the distance $d_G(u,v)$ between two vertices $u$ and $v$ of $G$ is the length of a shortest path between $u$ and $v$. If there is no ambiguity,  we will use the notation $d(u,v)$ instead of $d_G(u,v)$.  The diameter of $G$ is defined as $D(G)=\max_{u,v\in V}\{d(u,v)\}$.  Given $u,v\in V$,  $u\sim v$ means that $u$ and $v$ are adjacent vertices. Given a set of vertices $S=\{v_1,v_2,...,v_k\}$ of a connected graph $G$, the {\it metric representation} of a vertex $v\in V$ with respect to $S$ is the vector $r(v|S)=(d(v,v_1),d(v,v_2),...,d(v,v_k))$. We say that $S$ is a {\it resolving set} for $G$ if for every pair of distinct vertices $u,v\in V$, $r(u|S)\ne r(v|S)$. The {\it metric dimension} of $G$ is the minimum cardinality of any resolving set for $G$, and it is denoted by $dim(G)$.

Let $G$ and $H$ be two graphs of order $n_1$ and $n_2$, respectively. The corona product $G\odot H$ is defined as the graph obtained from $G$ and $H$ by taking one copy of $G$ and $n_1$ copies of $H$ and joining by an edge each vertex from the $i^{th}$-copy of
$H$ with the $i^{th}$-vertex of $G$. We will denote by $V=\{v_1,v_2,...,v_n\}$   the set of vertices of $G$ and by $H_i=(V_i,E_i)$   the copy of $H$ such that $v_i\sim v$ for every $v\in V_i$. Notice that the corona graph $K_1\odot H$ is isomorphic to the join graph $K_1+H$. For any integer  $k\ge 2$, we define the graph $G\odot^k H$
recursively from $G\odot H$ as $G\odot^k H=(G\odot^{k-1} H)\odot H$. We also note that the order of $G\odot^k H$ is $n_1(n_2+1)^{k}$.

\section{Metric dimension of corona product graphs}

We begin by presenting the following useful facts.

\begin{lemma}\label{first-remark}
Let $G=(V,E)$ be a connected graph of order $n\ge 2$ and let $H$ be a graph of order at least two. Let $H_i=(V_i,E_i)$ be the subgraph of $G\odot H$ corresponding to the $i^{th}$-copy of $H$.

\begin{itemize}

\item[{\rm (i)}]  If $u,v\in V_i$, then $d_{G\odot H}(u,x)=d_{G\odot H}(v,x)$ for every vertex $x$ of $G\odot H$ not belonging to $V_i$.


\item[{\rm (ii)}] If  $S$ is a resolving set for $G\odot H$, then  $V_i\cap S\neq \emptyset$ for every $i  \in \{1,...,n\}$.

\item[{\rm (iii)}]  If $S$ is a resolving set for $G\odot H$ of minimum cardinality, then  $V\cap S=\emptyset$.

\item[{\rm (iv)}]  If $H$ is a connected graph and $S$ is a resolving set for $G\odot H$, then for every $i\in \{1,..,n\}$, $S\cap V_i$ is a resolving set for $H_i$.
\end{itemize}

\end{lemma}

\begin{proof} (i) Let $y=v_i\in V$. The result directly follows from the fact that $d_{G\odot H}(u,x)=d_{G\odot H}(u,y)+d_{G\odot H}(y,x)=d_{G\odot H}(v,y)+d_{G\odot H}(y,x)=d_{G\odot H}(v,x)$.


(ii) We suppose $V_i\cap S= \emptyset$ for some $i\in \{1,...,n\}$. Let $x,y\in V_i$. By (i) we have $d_{G\odot H}(x,u)=d_{G\odot H}(y,u)$ for every vertex $u\in S$, which is a contradiction.

(iii) We will show that $S'=S-V$ is a resolving set for $G\odot H$.  Now let $x,y$ be two different vertices of  $G\odot H$.
We have the following cases.

Case 1: $x,y\in V_i$. By (i) we conclude that there exist $v\in V_i\cap S'$ such that  $d_{G\odot H}(x,v)\ne d_{G\odot H}(y,v)$.

Case 2: $x\in V_i$ and $y\in V_j$, $i\ne j$. Let $v\in V_i\cap S'$. Then we have $d_{G\odot H}(x,v)\le 2<3\le d_{G\odot H}(y,v)$.

Case 3: $x,y\in V$. Let $x=v_i$ and let $v\in V_i\cap S'$. Then  we have $d_{G\odot H}(x,v)=1<1+d_{G\odot H}(y,x)=d_{G\odot H}(y,v)$.

Case 4: $x\in V_i$ and $y\in V$. If $x\sim y$, then $y=v_i$. Let $v_j\in V$, $j\ne i$, and let $v\in V_j\cap S'$. Then we have $d_{G\odot H}(x,v)=1+d_{G\odot H}(y,v)>d_{G\odot H}(y,v)$.
For $x\not\sim y=v_l$ we take $v\in V_l\cap S'$ and we obtain $d_{G\odot H}(x,v)=d_{G\odot H}(x,y)+d_{G\odot H}(y,v)>d_{G\odot H}(y,v)$.

Therefore, $S'$ is a resolving set for $G\odot H$.

(iv) Let $S_i=S\cap V_i$. For $x\in S_i$ or $y\in S_i$  the result is straightforward.  We suppose $x,y\in V_i-S_i$. Since $S$ is a resolving set for $G\odot H$, we have $r(x|S)\ne r(y|S)$. By (i), $d_{G\odot H}(x,u)=d_{G\odot H}(y,u)$ for every vertex $u$ of $G\odot H$ not belonging to $V_i$. So, there exists $v\in S_i$ such that $d_{G\odot H}(x,v)\ne d_{G\odot H}(y,v)$. Thus, either ($v\sim x$ and $v\not\sim y$) or ($v\not\sim x$ and $v \sim y$). In the first case we have $d_{G\odot H}(x,v)=d_{H_i}(x,v)=1$ and $d_{G\odot H}(y,v)=2\le d_{H_i}(y,v)$. The case $v\not\sim x$ and $v\sim y$ is analogous. Therefore, $S_i$ is a resolving set for $H_i$.
\end{proof}

\begin{theorem}\label{general-lowerbound-dim}
Let $G$ and $H$ be two connected graphs of order $n_1\ge 2$ and  $n_2\ge 2$, respectively.
Then,  $$dim(G\odot^k H)\ge n_1(n_2+1)^{k-1}dim(H).$$
\end{theorem}

\begin{proof}
Let $S$ be a resolving set of minimum cardinality in $G\odot H$. From Lemma \ref{first-remark} (iii) we have that $S\cap V=\emptyset$. Moreover, by Lemma \ref{first-remark} (ii) we have that for every $i\in \{1,...,n_1\}$ there exist a nonempty set $S_i\subset V_i$ such that $S=\bigcup_{i=1}^{n_1}S_i$. Now, by using Lemma \ref{first-remark} (iv) we have that $S_i$ is a resolving set for $H_i$. Hence,  $dim(G\odot H)=|S|=\sum_{i=1}^{n_1}|S_i|\ge \sum_{i=1}^{n_1}dim(H)=n_1dim(H)$. As a result, the lower bound follows.
\end{proof}

\begin{theorem}\label{theo-corona-k}
Let $G$ be a connected graph of order $n_1\ge 2$ and let $H$ be a graph of order $n_2\ge 2$. If $D(H)\le 2$, then
$$dim(G\odot^k H)=n_1(n_2+1)^{k-1}dim(H).$$
\end{theorem}

\begin{proof}
 Let $S_i\subset V_i$ be a resolving set for $H_i$ and let $S=\bigcup_{i=1}^{n_1}S_i$. We will show that $S$ is a resolving set for $G\odot H$. Let us consider two different vertices $x,y$ of $G\odot H$. We have the following cases.

Case 1: $x,y\in V_i$. Since $D(H_i)\le 2$, we have that $r(x|S_i)\ne r(y|S_i)$ leads to $r(x|S)\ne r(y|S)$.

Case 2: $x\in V_i$ and $y\in V_j$, $i\ne j$. Let $v\in S_i$. Hence we have $d(x,v)\le 2<3\le d(y,v)$.

Case 3: $x,y\in V$. Let  $x=v_i$. Then for every vertex $v\in S_i$ we have $d(x,v)=1<d(y,x)+1=d(y,v)$.

Case 4: $x\in V_i$ and $y\in V$. If $x\sim y$, then let $v\in S_j$, for some $j\ne i$. So  we have $d(x,v)=1+d(y,v)>d(y,v)$. Moreover, if $x\not\sim y=v_j$, for $v\in S_j$  we have $d(x,v)=d(x,y)+d(y,v)>d(y,v)$.

Thus, for every different vertices $x,y$ of $G\odot H$, we have $r(x|S)\ne r(y|S)$, as a consequence, $dim(G\odot H)\le n_1dim(H)$. Therefore, we have $dim(G\odot^k H)\le n_1(n_2+1)^{k-1}dim(H)$.
By Theorem \ref{general-lowerbound-dim} we conclude the proof.
\end{proof}

In order to show a consequence of the above theorem we present the following well known result, where $K_t$ denotes a complete graph of order $t$,  $K_{s,t}$ denotes a complete bipartite graph of order $s+t$ and $N_t$ denotes an empty graph of order $t$.

\begin{lemma}{\em\cite{chartrand}}
Let $G$ be a connected graph of order $n\ge 4$. Then $dim(G)=n-2$ if and only if $G = K_{s,t}$,  $(s,t\ge 1)$, $G = K_s + N_t$, $(s\ge 1$, $t\ge 2)$, or $G = K_s + (K_1 \cup K_t)$, $(s, t\ge 1)$.
\end{lemma}

\begin{corollary}\label{dim-n-2}
Let $G$ be a connected graph of order $n_1\ge 2$ and let $H$ be a graph of order $n_2\ge 4$ and diameter $D(H)\le 2$. Then $$dim(G\odot^k H)=n_1(n_2+1)^{k-1}(n_2-2)$$ if and only if $H = K_{s,t}$,  $(s,t\ge 1)$; $H = K_s + N_t$, $(s\ge 1$, $t\ge 2)$, or $H = K_s + (K_1 \cup K_t)$, $(s, t\ge 1)$.
\end{corollary}

We recall that the wheel graph of order $n+1$ is defined as $W_{1,n}=K_1\odot C_{n}$, where $K_1$ is the singleton graph and $C_{n}$ is the cycle graph of order $n$. The metric dimension of the wheel $W_{1,n}$ was obtained by Buczkowski et. al. in \cite{buczkowski}.

\begin{remark}{\em \cite{buczkowski}}\label{dim-wheel} Let $W_{1,n}$ be a wheel graph. Then
$$dim(W_{1,n})=\left\{\begin{array}{ll}
3 & \textrm{for $n=3,6$,}\\
2 & \textrm{for $n=4,5$,}\\
{\left\lfloor{\frac{2n+2}{5}}\right\rfloor} & \textrm{otherwise.}\\
\end{array}\right.$$
\end{remark}

The fan graph $F_{n_1,n_2}$ is defined as the graph join $N_{n_1}+P_{n_2}$, where $N_{n_1}$ is the empty graph of order $n_1$ and $P_{n_2}$ is the path graph of order $n_2$. The case $n_1=1$ corresponds to the usual fan graphs. Notice that, for the metric dimension of fan graphs, it is possible to find an equivalent result to Remark \ref{dim-wheel} which was obtained by Caceres et. al.  in \cite{pelayo2}.

\begin{remark}\label{dim-fan}{\em \cite{pelayo2}}
Let $F_{1,n}$ be a fan graph. Then
$$dim(F_{1,n})=\left\{\begin{array}{ll}
1 & \textrm{for $n=1$,}\\
2 & \textrm{for $n=2,3$,}\\
3 & \textrm{for $n=6$,}\\
{\left\lfloor{\frac{2n+2}{5}}\right\rfloor} & \textrm{otherwise.}\\
\end{array}\right.
$$
\end{remark}

As a particular case of the Theorem \ref{theo-corona-k} we obtain the following results.

\begin{corollary} Let $G$ be a connected graph of order $n_1\ge 2$. If $H$ is a wheel graph or a fan graph of order $n_2\ge 8$, then
$$dim(G\odot^k H)=n_1(n_2+1)^{k-1}\left\lfloor{\frac{2n_2}{5}}\right\rfloor.$$
\end{corollary}

\begin{theorem}\label{general-bound-dim}
Let $G$ be a connected graph of order $n_1\ge 2$ and let $H$ be a graph of order $n_2\ge 2$. Let  $\alpha$ be the number of connected components of $H$ of order greater than one and let $\beta$ be the number of isolated vertices of $H$. Then
$$dim(G\odot^k H)\le \left\{\begin{array}{ll}
n_1(n_2+1)^{k-1}(n_2-\alpha-1) & \textrm{for  $\alpha\ge 1$ and $\beta\ge 1$,}\\
\\
n_1(n_2+1)^{k-1}(n_2-\alpha) & \textrm{for  $\alpha\ge 1$ and $\beta = 0$,}\\
\\
n_1(n_2+1)^{k-1}(n_2-1) & \textrm{for  $\alpha =0$.}
\end{array}\right.
$$
\end{theorem}

\begin{proof}

We suppose $\alpha\ge 1$ and $\beta\ge  1$. Let $A_i$ be the set of vertices of $G\odot H$ formed by all but one of the vertices per each of the $\alpha$ connected components of $H_i$. If $\beta\ge 2$ we define  $B_i$ to be the set  of vertices of $G\odot H$ formed by all but one of the isolated vertices of $H_i$.  If $\beta=1$ we assume $B_i=\emptyset$. Let us show that $S=\cup_{j=1}^{n_1} (A_j\cup B_j)$ is a resolving set for $G\odot H$. Let $x,y$ be two different vertices of $G\odot H$. We suppose $x,y\notin S$.  We have the following cases.

Case 1. $x=v_i\in V$ and  $y\in V_i$.  For every vertex $u\in V_j\cap S$, $j\ne i$, we obtain  $d(y,u)=d(y,x)+d(x,u)>d(x,u)$.

case 2.  $x=v_i\in V$ and $y\not\in V_i$.  For every $v\in S\cap V_i$  we have $d(x,v)=1< d(y,v)$.

Case 3.  $x\in V_i$ and $y\in V_j$, $j\ne i$.  For every $u\in V_i\cap S$ we have $d(x,u)\le 2<3\le d(y,u)$.

Case 4.  $x,y\in V_i$. We consider, without loss of generality, that  $x$  is not an isolated vertex in $H_i$. Then there exists $v\in V_i\cap S$ such that  $v\sim x$, so  $d(x,v)=1< 2=d(y,v)$.

Thus, for every two different vertices $x,y$ of $G\odot H$, we obtain $r(x|S)\ne r(y|S)$ and, as a consequence, $dim(G\odot H)\le n_1(n_2-\alpha-1).$

As above, if   $\beta =0$  then we take $S=\cup_{j=1}^{n_1}  A_j$ and we obtain $dim(G\odot H)\le n_1(n_2-\alpha)$  and if $\alpha =0$, then we take $S=\cup_{j=1}^{n_1}  B_j$ and we obtain $dim(G\odot H)\le n_1(n_2-1).$ Note that if $\alpha =0$, then it is not necessary to consider Case 4.  Thus, the result  follows.
\end{proof}

\begin{corollary}\label{dim-n-1-N_n2}
Let $G$ be a connected graphs of order $n_1\ge 2$ and let $H$ be an unconnected graph of order $n_2\ge 2$. Then $$dim(G\odot^k H)=n_1(n_2+1)^{k-1}(n_2-1)$$ if and only if $H\cong N_{n_2}$.
\end{corollary}
\begin{proof}
In \cite{baskoro} the authors showed that  $dim(G\odot N_{n_2})=n_1(n_2-1)$. Hence, $dim(G\odot^k N_{n_2})=n_1(n_2+1)^{k-1}(n_2-1)$.
Moreover, by the above theorem, if $H$ is unconnected and $H \not\cong N_{n_2}$, then $dim(G\odot^k H)\le n_1(n_2+1)^{k-1}(n_2-2)$.
\end{proof}

\begin{theorem}\label{dim-n-1-kn}
Let $G$ and $H$ be two connected graphs of order $n_1\ge 2$ and $n_2\ge 3$, respectively. Then $$dim(G\odot^k H)=n_1(n_2+1)^{k-1}(n_2-1)$$ if and only if $H\cong K_{n_2}$. Moreover, if $H\not\cong K_{n_2}$, then
$$dim(G\odot^k H)\le n_1(n_2+1)^{k-1}(n_2-2).$$
\end{theorem}

\begin{proof}
Since $dim(K_{n_2})=n_2-1$, by  Theorem \ref{theo-corona-k} we conclude $dim(G\odot^k K_{n_2})=n_1(n_2+1)^{k-1}(n_2-1)$. On the contrary, we suppose $H\not\cong K_{n_2}$. Given a set $X$ of vertices of $H$ and a vertex $v$ of $H$, $N_X(v)$ denotes the set of neighbors that $v$ has in $X$:
$N_X(v)=\{u\in X:\;u\sim v\}$. Given two vertices $a,b$ of $H$, let $X_{a,b}$ be the set formed by all vertices of $H$ different from $a$ and $b$. Since $H$ is a connected graph and $H\ne K_{n_2}$, there exist at least two vertices $a,b$ of $H$ such that $N_{X_{a,b}}(a)\ne N_{X_{a,b}}(b)$. Let $a_i, b_i$ be the vertices corresponding to $a,b$, respectively, in the $i^{th}$-copy $H_i=(V_i,E_i)$ of $H$. Let $S=\cup_{i=1}^{n_2}(V_i-\{a_i,b_i\})$. We will show that $S$ is a resolving set for $G\odot H$. Let $x,y$ be two different vertices of $G\odot H$ such that $x,y\not\in S$. We have the following cases.

 Case 1. $x=a_i$ and $y=b_i$. Since $N_{X_{a,b}}(a)\ne N_{X_{a,b}}(b)$ we have $r(x|S)\ne r(y|S)$.

 Case 2.  $x=v_i\in V$ and $y\in V_i$.
 For every $v\in V_j-\{a_j,b_j\}$, $j\ne i$,  we have $d(y,v)=d(y,x)+d(x,v)>d(x,v)$. If $x\in V_i$ and $y\in V_j$, $j\ne i$, then for every $v\in V_i-\{a_i,b_i\}$ we have $d(x,v)\le 2<3\le d(y,v)$.

 Case 3. $x,y\in V$. Say $x=v_i$.   Then for every $v\in  V_i-\{a_i,b_i\}$ we have  $d(x,v)=1< d(y,v)$.

 Hence, for every two different vertices $x,y$ of $G\odot H$, we obtain $r(x|S)\ne r(y|S)$. Thus, $dim(G\odot H)\le n_1(n_2-2).$
 Therefore, the result follows.
\end{proof}

As we have shown in Corollary \ref{dim-n-2}, the above bound is tight.

\begin{theorem}\label{supCoronaK1} Let $G$ be a connected graph of order $n_1\ge 2$ and let $H$ be a graph of order $n_2\ge 2$. Then
$$dim(G\odot^k H)\le n_1(n_2+1)^{k-1}dim(K_1\odot H).$$
\end{theorem}
\begin{proof}
 We denote by $K_1\odot H_i$  the subgraph of $G\odot H$, obtained by joining the vertex $v_i\in V$ with all vertices of  $H_i$. For every $v_i\in V$, let $B_i$ be a resolving set of minimum cardinality of $K_1\odot H_i$ and let $B=\bigcup_{i=1}^{n_1}B_i$. By Lemma \ref{first-remark} (iii) we have that $v_i$ does not belong to any resolving set of minimum cardinality for $K_1\odot H_i$. So, $B$ does not contain any vertex from $G$.  We will show that $B$ is a resolving set for $G\odot H$. Let $x,y$ be two different vertices in $G\odot H$. We consider the following cases.

Case 1: $x,y\in V_i$. There exists $u\in B_i$ such that $d_{K_1\odot H_i}(x,u)\ne d_{K_1\odot H_i}(y,u)$, which leads to $d_{G\odot H}(x,u)\ne d_{G\odot H}(y,u)$.

Case 2: $x\in V_i$ and $y\in V_j$, $i\ne j$. Let $v\in B_i$. We have $d_{G\odot H}(x,v)\le 2<3\le d_{G\odot H}(y,v)$.

Case 3: $x, y\in V$. Suppose now that $x$ is adjacent to the vertices of $H_i$. Hence, for every vertex $v\in B_i$ we have $d_{G\odot H}(x,v)=1<d_{G\odot H}(y,x)+1=d_{G\odot H}(y,v)$.

Case 4: $x\in V_i$ and $y\in V$. If $x\sim y$, then for every vertex $v\in B_j$, with $j\ne i$, we have $d_{G\odot H}(x,v)=1+d_{G\odot H}(y,v)>d_{G\odot H}(y,v)$. Now, let us assume that $x\not\sim y$. Hence, there exists $v\in B_j$ adjacent to $y$, with $j\ne i$. So, we have $d_{G\odot H}(x,v)=d_{G\odot H}(x,y)+1=d_{G\odot H}(x,y)+d_{G\odot H}(y,v)>d_{G\odot H}(y,v)$.

Thus, for every two different vertices $x,y$ of $G\odot H,$ we have $r(x|S)\ne r(y|S)$ and, as a consequence, $dim(G\odot H)\le n_1dim(K_1\odot H)$.
Therefore, the result follows.
\end{proof}

\begin{theorem}\label{exact-values-dim}Let $G$ be a connected graph of order $n_1\ge 2$ and let $H$ be a graph of order $n_2\ge 7$. If $D(H)\ge 6$ or $H$ is a cycle graph, then
$$dim(G\odot^k H)=n_1(n_2+1)^{k-1}dim(K_1\odot H).$$
\end{theorem}

\begin{proof}
Let $S$ be a resolving set of minimum cardinality in $G\odot H$. By Lemma \ref{first-remark} (iii)  we have $S\cap V=\emptyset$, as a consequence, $S=\cup_{i=1}^{n_1}S_i$,
where $ S_i\subset V_i$. Notice that, by Lemma \ref{first-remark} (ii),
$S_i\ne \emptyset$ for every $i \in \{1,...,n_1\}$.
 Now we differentiate two cases in order to show that  $r(x|S_i)\ne (1,...,1)$ for every  $x\in V_i-S_i$.

Case 1.  $H$ is a cycle graph of order $n_2\ge 7$.   If $r(a|S_i)=(1,1)$ for some $a\in V_i-S_i$, then, since $n_2\ge 7$,  there exist two vertices $x,y\in V_i-S_i$ such that $d_{H_i}(x,v)>1$ and $d_{H_i}(y,v)>1$, for every $v\in S_i$. Hence, $d_{G\odot H}(x,v)=d_{G\odot H}(y,v)=2$ for every $v\in S_i$, which is a contradiction because, by Lemma \ref{first-remark} (i),   $d_{G\odot H}(x,v)=d_{G\odot H}(y,v)$ for every vertex $u$ of $S$ not belonging to $S_i$.

Case 2. $D(H)\ge 6$.  Let $x,y\in V_i-S_i$. Since $S$ is a resolving set for $G\odot H$, we have $r(x|S)\ne r(y|S)$. As we have noted before, by Lemma \ref{first-remark} (i) we have that $d_{G\odot H}(x,u)=d_{G\odot H}(y,u)$ for every vertex $u$ of $G\odot H$ not belonging to $V_i$. So, there exists $v\in S_i$ such that $d_{G\odot H}(x,v)\ne d_{G\odot H}(y,v)$ and,  as a consequence,
either ($v\sim x$ and $v\not\sim y$) or ($v\not\sim x$ and $v \sim y$). Now we suppose that there exists a vertex  $a\in V_i-S_i$ such that $r(a|S_i)=(1,1,...1)$. If there exists a vertex $b\in V_i-S_i$ such that $d_{H_i}(b,u)>1$, for every $u\in S_i$, then for every $w\in V_i-(S_i\cup \{a,b\})$, there exists $v\in S_i$ such that $w\sim v$. Then  $D(H_i)\le 5$. Moreover, if for every  $b\in V_i-S_i$  there exists $v_b\in S_i$ such that $v_b\sim b$, then $D(H)\le 4$. Therefore, if $D(H)\ge 6$, then  $r(a|S_i)\ne(1,1,...1)$ for every $a\in V_i-S_i$.

Now, we denote by $K_1\odot H_i$  the subgraph of $G\odot H$, obtained by joining the vertex $v_i\in V$ with all vertices of the $i^{th}$-copy of $H$. In both the above cases
we have $r(v_i|S_i)=(1,1,...,1)\ne r(x|S_i)$ for every $x\in V_i-S_i$, so  $S_i$ is a resolving set for $K_1\odot H_i$. Hence, $dim(K_1\odot H_i)\le |S_i|$, for every $i\in \{1,...,n_1\}$.  Thus, $dim(G\odot H)\ge n_1dim(K_1\odot H_i)$ and, as a consequence,  $dim(G\odot^k H)\ge n_1(n_2+1)^{k-1}dim(K_1\odot H)$. We conclude the proof by Theorem \ref{supCoronaK1}.
\end{proof}

\begin{corollary}Let $G$ be a connected graph of order $n_1\ge 2$.
\begin{itemize}
\item[{\rm (i)}]If $n_2\ge 7$, then $dim(G\odot^k C_{n_2})=n_1(n_2+1)^{k-1}\left\lfloor{\frac{2n_2+2}{5}}\right\rfloor.$
\item[{\rm (ii)}]If $n_2\ge 7$, then $dim(G\odot^k P_{n_2})=n_1(n_2+1)^{k-1}\left\lfloor{\frac{2n_2+2}{5}}\right\rfloor.$
\end{itemize}
\end{corollary}

All our previous results concern to $G\odot H$ for $H$ of order at least two. Now we consider the case $H\cong  K_1$.  We obtain a general bound for  $dim(G\odot^k K_1)$ and, when $G$ is a tree,  we give the exact value for this parameter.

\begin{claim}\label{neib-dif-dist}
Let $G$ be a simple graph. If $v$ is a vertex of degree greater than one in $G$, then for every vertex $u$ adjacent to $v$ there exists a vertex $x\ne u,v$ of $G$, such that $d(v,x)\ne d(u,x)+1$.
\end{claim}

The following lemma obtained in \cite{buczkowski} is useful to obtain the next result.

\begin{lemma}{\em\cite{buczkowski}}\label{one-pendant-edge}
If $G_1$ is a graph obtained by adding a pendant edge to a nontrivial connected graph $G$, then $dim(G)\le dim(G_1)\le dim(G) + 1.$
\end{lemma}

\begin{theorem}
For every connected graph $G$ of order $n\ge 2$,  $$dim(G\odot^k K_1)\le 2^{k-1}n-1.$$
\end{theorem}

\begin{proof}
If $G\cong K_2$, then $dim(K_2\odot K_1)=dim(P_4)=1$. So, let us suppose $G\not\cong  K_2$. Let us suppose, without loss of generality, that $v_n$ is a vertex of degree greater than one in $G$ and let $S=V-\{v_n\}$. For every $i\in \{1,...,n\}$, let $u_i$ be the pendant vertex of $v_i$ in $G\odot K_1$. We will show that $S$ is a resolving set for $G\odot K_1$. Let $x,y$ be two different vertices of $G\odot K_1$. If $x=u_i$ and $y=u_j$, $i\ne j$, then we have either $i\ne n$ or $j\ne n$. Let us suppose for instance $i\ne n$. So, we obtain that $d(x,v_i)=1\ne d(y,v_i)$. On the other hand, if $x=v_n$ and $y=u_i$, then let us suppose $d(x,v_i)=1$. Since $v_n$ is a vertex of degree greater than one in $G$, by Claim \ref{neib-dif-dist}, there exists a vertex $v_j\in S$ such that $d(x,v_j)\ne d(v_i,v_j)+1$. So, we have
$
d(x,v_j)\ne d(v_i,v_j)+1
= d(v_i,v_j)+d(u_i,v_i)
=d(y,v_i)+d(v_i,v_j)
=d(y,v_j).
$
Therefore, for every different vertices $x,y$ of $G\odot K_1$ we have $r(x|S)\ne r(y|S)$ and, as a consequence, $dim(G\odot K_1)\le n-1$. Therefore, $dim(G\odot^k K_1)\le 2^{k-1}n-1$.
\end{proof}

By Lemma \ref{one-pendant-edge} we have $dim(K_n\odot K_1)\ge dim(K_n)=n-1$. Thus, for $k=1$ the above bound is achieved for the graph $G=K_n$.

 To present the next result, we need additional definitions. A vertex of degree at least $3$ in a graph $G$ will be called a {\it major vertex} of $G$.
Any vertex $u$ of degree one is said to be a {\it terminal vertex} of a major vertex $v$  if
$d(u, v)<d(u,w)$ for every other major vertex $w$ of $G$. The {\it terminal degree}  of
a major vertex $v$ is the number of terminal vertices of $v$. A major vertex $v$ is an
{\it exterior major vertex}  if it has positive terminal degree. Given a graph $G$,  $n_1(G)$ denotes the number of vertices of degree one and  $ex(G)$ denotes the number
of exterior major vertices of $G$.

\begin{lemma}{\em \cite{chartrand,harary,leaves-trees}}\label{value-dim-trees}
If $T$ is a tree that is not a path, then
$dim(T) = n_1(T) - ex(T).$
\end{lemma}

\begin{theorem}
For any tree $T$  of order $n\ge 3$,
$$dim(T\odot^k K_1)=\left\{\begin{array}{ll}
n_1(T) & \textrm{for $k=1$,}\\
\\
2^{k-2}n & \textrm{for $k\ge 2$.}
\end{array}\right.
$$
\end{theorem}

\begin{proof}
If $T$ is a path of order $n\ge 3$, then we have $dim(T\odot K_1)=2=n_1(T)$. Now, if $T$ is not a path, then by using Lemma \ref{value-dim-trees}, since $T\odot K_1$ is a tree, $n_1(T\odot K_1)=n$ and $ex(T\odot K_1)=n-n_1(T)$, we obtain the result for $k=1$. Since for every tree $T$ of order $n$ we have $n_1(T\odot K_1)=n$, we obtain the result for $k\ge 2$.
\end{proof}

\section*{Acknowledgements}
The research was partially done while the first author was at
Gda\'nsk University of Technology, Poland,  supported by \lq\lq Fundaci\'o Ferran Sunyer i Balaguer'', Catalunya, Spain. This work was partly supported  by the Spanish Ministry of Education through projects TSI2007-65406-C03-01 \lq\lq E-AEGIS" and CONSOLIDER INGENIO 2010 CSD2007-00004 \lq\lq ARES''.

\end{document}